\documentclass{amsart}
\usepackage{amsmath,amssymb,amsthm}

\numberwithin{equation}{section}
\newtheorem{theorem}{Theorem}[section]
\newtheorem{lemma}[theorem]{Lemma}

\newtheorem{corollary}[theorem]{Corollary}
\theoremstyle{definition}
\newtheorem{example}[theorem]{Example}
\theoremstyle{remark}
\newtheorem{remark}[theorem]{Remark}

\newcommand{\R}{\mathbb R}
\newcommand{\C}{\mathbb C}
\newcommand{\E}{\mathbb E}
\newcommand{\Ran}{\operatorname{Ran}}
\newcommand{\Sym}{\operatorname{Sym}}
\newcommand{\spec}{\sigma}
\newcommand{\cC}{\mathcal C}
\newcommand{\cL}{\mathcal L}

\title[Mean-square stability with delayed diffusion]{Mean-square stability of multidimensional stochastic differential equations with delayed diffusion}
\author{John A.D. Appleby}
\address{School of Mathematical Sciences, Dublin City University, Glasnevin, Dublin 9, Ireland}
\email{john.appleby@dcu.ie}
\author{Wilson Tai}
\address{School of Mathematical Sciences, Dublin City University, Glasnevin, Dublin 9, Ireland}
\email{wilson.tai2@mail.dcu.ie}
\date{}

\begin{document}
\begin{abstract}
We give a spectral characterisation of mean-square asymptotic stability for a multidimensional stochastic delay differential equation with instantaneous drift and delayed diffusion. A variation-of-constants formula reduces the stochastic equation to an exact deterministic matrix Volterra equation for the second moments. Its kernel has an explicit Laplace transform, and the associated characteristic matrix is related, on the symmetric tensor space, to the characteristic matrix of a classical retarded differential equation. This relation yields a dominant asymptotic expansion of the second-moment matrix, including polynomial factors arising from higher-order resolvent poles and explicit dependence on the initial history. We identify an accessibility condition under which a dominant characteristic mode is excited, and show that simultaneous cancellation of all dominant modes is non-generic. Under this condition we obtain a necessary and sufficient criterion for mean-square asymptotic stability. The sufficient stability conclusion also holds for independent random initial histories with finite second moments. The accessibility condition is automatic when the delayed diffusion matrix is nonsingular.
\end{abstract}
\maketitle

\section{Introduction}
Stability questions for deterministic delay equations have long been tied to the location of the roots of their characteristic equations; see, for example, Bellman and Cooke \cite{BellmanCooke}, Hale and Verduyn Lunel \cite{Hale}, Diekmann et al.\ \cite{Diekmann}, and the monographs of Kolmanovskii and co-authors \cite{KolmanovskiiNosov,KolmanovskiiMyshkis}. For stochastic differential equations, Lyapunov functionals, Razumikhin arguments, and linear matrix inequalities provide powerful stability criteria. In many settings, however, these criteria are sufficient conditions rather than exact spectral tests; see Mao \cite{MaoRaz,MaoBook}, Mao and Shah \cite{MaoShah}, Shaikhet \cite{Shaikhet}, and Bao, Yin and Yuan \cite{Bao}.

For linear stochastic systems without delay, second moments can be lifted to a deterministic linear system on a tensor space. This is the familiar Kronecker-product approach developed in the stability analysis of random linear systems; see Bellman \cite{Bellman54,Bellman59}, Nevelson and Khasminskii \cite{NevelsonKhasminskii}, and Khasminskii \cite{Khasminskii}. With delay, the corresponding lifted dynamics are infinite-dimensional, and the relationship between their characteristic roots and mean-square stability requires care.

Several papers are close in spirit to the present work, but they address different levels of the mean-square problem. Appleby, Mao and Riedle \cite{AMR} obtain a necessary and sufficient characterisation of the mean-square asymptotic behaviour of a scalar stochastic functional differential equation with finite memory; see also Mao and Riedle \cite{MaoRiedle} for stochastic Volterra integro-differential equations and Appleby and Lawless \cite{ApplebyLawless} for perturbed linear stochastic functional differential equations. Their criteria and rates are expressed in terms of the differential resolvent of the underlying deterministic equation. This is a general characterisation, but the resolvent need not be available explicitly. Lei and Mackey \cite{LeiMackey} use Laplace-transform methods to obtain necessary and sufficient conditions for moment boundedness in a scalar stochastic delay equation, with an application to a hematopoietic stem-cell model. Mackey and Nechaeva \cite{MackeyNechaeva} study first- and second-moment stability, including changes in stability and oscillatory behaviour under white and coloured noise.

For multidimensional equations, Wang, Li and Lei \cite{WangLiLei} develop a Laplace-transform framework for the second-moment boundedness problem with one discrete delay. Their general \(N\)-dimensional construction produces a characteristic function involving transforms of products of entries of the deterministic fundamental matrix. Their explicit analysis is then specialised to a two-dimensional system with decoupled stochastic terms. In that setting their Theorem~4.3 gives boundedness and exponential convergence of the second-moment matrix to a constant when all characteristic roots lie in the open left half-plane, and gives an unboundedness conclusion under a positive-root hypothesis and an additional non-degeneracy assumption. The critical spectral case is left open. Their Theorem~4.7 supplies a practical sufficient condition for boundedness. They subsequently treat distributed-delay equations in \cite{WangDistributed}, again deriving characteristic functions and sufficient conditions for second-moment boundedness or unboundedness. These papers are important predecessors for the use of Laplace transforms and characteristic functions. They do obtain boundedness conclusions and, in the stable case of the two-dimensional discrete-delay paper, convergence to a limiting second-moment matrix. They do not, however, give the dominant polynomial-exponential expansion of the second-moment matrix, or a mode-by-mode description of the leading coefficients and their dependence on the initial history.

A recent complementary direction is the computation of stability boundaries: Iklodi and Dankowicz \cite{IklodiDankowicz} derive algebraic conditions for second-moment stability boundaries of linear time-invariant stochastic delay equations. Their objective is to locate the boundary efficiently, rather than to obtain a time-domain dominant asymptotic expansion of the second-moment matrix.

Buckwar and Sickenberger \cite{BuckwarSickenberger} develop vectorisation and Kronecker-product methods for multidimensional linear stochastic differential systems without delay. Buckwar and Notarangelo \cite{BuckwarNotarangelo} extend this line of work to stochastic delay differential systems and state necessary and sufficient mean-square stability conditions for the class they consider. Thus, neither tensor lifting nor the general goal of obtaining an exact stability test for multidimensional stochastic delay equations is new. There is, however, an important qualification concerning the covariance-based spectral approach in \cite{MackeyNechaeva,BuckwarNotarangelo}. Appleby \cite[Section~4.4]{ApplebyVolterra} examines the separable exponential covariance ansatz used in these works and gives a scalar delayed-diffusion example for which its prediction is inconsistent with the exact covariance asymptotics. We discuss this point below; it is a reason to avoid relying on that ansatz, rather than a basis for asserting that every stability conclusion in those papers is false.

The contribution here is more specific. For the model below, in which the drift is instantaneous and the delay occurs in the diffusion coefficient, we derive the second-moment Volterra equation directly from variation of constants and It\^o's isometry. Its kernel has an explicit Laplace transform, and its physical characteristic determinant on the symmetric tensor space is related to the characteristic matrix of a concrete retarded equation. This gives more than a stability test: it yields a dominant asymptotic expansion, including polynomial factors associated with higher-order resolvent poles, and identifies the initial-history coefficients that determine whether a dominant mode is actually present. We also prove a generic non-cancellation result and use accessibility to obtain a necessary and sufficient mean-square stability criterion. In contrast to the covariance-ansatz approach, the argument starts from an exact Volterra equation for the second moments and does not postulate a two-time exponential form for the covariance.

For clarity, the scalar equation
\[
dX(t)=aX(t)\,dt+\sigma X(t-\tau)\,dB(t)
\]
is a special case of our model. It is also the \(b=0\) case of the equation studied in \cite{ApplebyVolterra}. In that case Appleby proves that, for each fixed \(k\geq0\),
\[
\lim_{t\to\infty}\frac{\mathbb E[X(t)X(t+k)]}{\mathbb E[X(t)^2]}=e^{ak}.
\]
By contrast, the ansatz \(K(t,s)=\mathbb E[X(t)X(s)]\sim e^{\lambda t}e^{\lambda s}\) predicts that this ratio tends to \(e^{\lambda k}\), and hence forces \(\lambda=a\). The diagonal second-moment equation obtained from the same ansatz gives
\[
2\lambda=2a+\sigma^2e^{-2\lambda\tau},
\]
which is incompatible with \(\lambda=a\) whenever \(\sigma\ne0\). This concrete contradiction is the point of the caution in \cite[Section~4.4]{ApplebyVolterra}. It concerns the proposed separable covariance asymptotic in general; it should not be read as showing, without further analysis, that every stability result in \cite{MackeyNechaeva,BuckwarNotarangelo} fails.

We study the equation
\begin{equation}\label{eq:model}
	dX(t)=AX(t)\,dt+\bigl(\Sigma_1X(t)+\Sigma_2X(t-\tau)\bigr)\,dB(t),\qquad t>0,
\end{equation}
where $X(t)\in\R^d$, $A,\Sigma_1,\Sigma_2\in\R^{d\times d}$ are constant matrices, $\tau>0$, and $B$ is a one-dimensional Brownian motion. The initial history is a deterministic function $\psi\in C([-\tau,0];\R^d)$ (see Corollary~\ref{cor:random} for the extension to random initial histories). The restriction to instantaneous drift with constant coefficients and a single discrete delay is deliberate: it makes the deterministic fundamental matrix explicit while retaining a genuinely delayed stochastic feedback.

The main steps are as follows. First, variation of constants and It\^o's isometry reduce the second moments to a deterministic Volterra equation. Secondly, the Laplace transform of its kernel can be computed explicitly. The resulting characteristic matrix is related to that of a retarded equation on the tensor space. Thirdly, we restrict to the symmetric tensor subspace, which is the natural state space for second moments. Finally, we identify an accessibility condition for dominant modes and use it to prove an exact stability criterion.

It is useful to distinguish three spectral objects from the outset. The first is the spectrum of the instantaneous tensor generator $A_2=A\otimes I_d+I_d\otimes A$, which determines the initial half-plane of convergence of the Volterra-kernel transform. The second is the characteristic-root set of the lifted retarded equation restricted to the symmetric tensor space $S_d$. The third is the zero set of the physical characteristic determinant of the Volterra equation, also on $S_d$. The second and third objects coincide in the half-plane to the right of the spectral bound of $A_2$. The full tensor-space retarded equation is a useful algebraic representation, but its antisymmetric directions do not describe second moments and are not used in the stability criterion. We use the adjective ``physical'' to emphasise that the determinant is formed on the symmetric tensor space actually occupied by second moments.

\section{Preliminaries and reduction to a Volterra equation}
Let $\cC=C([-\tau,0];\R^d)$ with its supremum norm $\|\psi\|_{\cC}=\sup_{\theta\in[-\tau,0]}\|\psi(\theta)\|$, and fix $\psi\in\cC$. Standard existence and uniqueness results for linear stochastic functional differential equations ensure that \eqref{eq:model} has a unique adapted strong solution with finite second moments; see, for example, Mohammed \cite{Mohammed} and Mizel and Trutzer \cite{MizelTrutzer}.

Define the causal fundamental matrix
\[
r(t)=\begin{cases}e^{At},&t\geq0,\\0,&t<0,\end{cases}
\]
and set
\[
x_0(t)=\begin{cases}e^{At}\psi(0),&t\geq0,\\\psi(t),&-\tau\leq t<0.\end{cases}
\]
Thus $x_0(t)=e^{At}\psi(0)$ for nonnegative $t$. Write
\begin{equation}\label{eq:Ydef}
	Y(t)=\Sigma_1X(t)+\Sigma_2X(t-\tau),\qquad t\geq0,
\end{equation}
and define
\[
y_0(t)=\Sigma_1x_0(t)+\Sigma_2x_0(t-\tau),\qquad
G(t)=\Sigma_1r(t)+\Sigma_2r(t-\tau).
\]

\begin{theorem}\label{thm:variation}
	For $t\geq0$,
	\begin{align}
		X(t)&=x_0(t)+\int_0^t r(t-s)Y(s)\,dB(s),\label{eq:Xvariation}\\
		Y(t)&=y_0(t)+\int_0^t G(t-s)Y(s)\,dB(s).\label{eq:Yvariation}
	\end{align}
\end{theorem}
\begin{proof}
	The stochastic variation-of-constants formula for the equation with instantaneous drift gives
	\[X(t)=e^{At}\psi(0)+\int_0^t e^{A(t-s)}Y(s)\,dB(s),\]
	which is \eqref{eq:Xvariation}. Apply this identity at $t$ and $t-\tau$, using the convention $r(u)=0$ for $u<0$, and substitute into \eqref{eq:Ydef}. The stochastic integral involving $X(t-\tau)$ can then be written over $[0,t]$, giving \eqref{eq:Yvariation}.
\end{proof}

Put $A_2=A\otimes I_d+I_d\otimes A$ and
\[
R(t)=r(t)\otimes r(t)=e^{A_2t},\qquad t\geq0.
\]
For $\Re z>\alpha_A$, where
\begin{equation}\label{eq:alphaA}
	\alpha_A=\max\{\Re\lambda:\lambda\in\spec(A_2)\},
\end{equation}
we have
\begin{equation}\label{eq:Rhat}
	\widehat R(z)=(zI_{d^2}-A_2)^{-1}.
\end{equation}

The natural space for second moments is
\[
S_d=\operatorname{span}\{u\otimes u:u\in\R^d\}\subset\R^{d^2}.
\]
It is canonically identified with the space of symmetric $d\times d$ matrices and has dimension $d(d+1)/2$. If an operator $T$ on $\R^{d^2}$ leaves $S_d$ invariant, write $T_s$ for its restriction to $S_d$, and write $I_{S_d}$ for the identity on $S_d$. In particular, $K\otimes K$ leaves $S_d$ invariant for every $K\in\R^{d\times d}$ because $(K\otimes K)(u\otimes u)=(Ku)\otimes(Ku)\in S_d$. When taking Laplace transforms or considering eigenvalues, resolvents, and characteristic determinants, we extend $S_d$ and operators $T_s$ by $\C$-linearity to the complexification $S_d^{\C}=S_d+iS_d\subset\C^{d^2}$, retaining the notation $S_d$ and $T_s$ for the complexified space and operators.

We use the following standard matrix-valued Paley--Wiener criterion for Volterra resolvents. The finite-dimensional matrix formulation follows from the matrix-valued Wiener/Paley--Wiener theorem; see Gripenberg, Londen and Staffans \cite[Chapter~5, in particular the resolvent and inversion results for convolution equations]{Gripenberg}. The statement below is the precise version needed here.

\begin{lemma}[Paley--Wiener criterion]\label{lem:PW}
	Let $K\in L^1(0,\infty;\cL(S_d))$, and let $\rho$ be its Volterra resolvent, defined by $\rho-K*\rho=K$. Then $\rho\in L^1(0,\infty;\cL(S_d))$ if and only if
	\[
	\det\bigl(I_{S_d}-\widehat K(z)\bigr)\ne0\qquad(\Re z\geq0).
	\]
\end{lemma}
\begin{proof}
	Taking Laplace transforms in the resolvent equation gives
	\[\widehat\rho(z)=\bigl(I_{S_d}-\widehat K(z)\bigr)^{-1}\widehat K(z)\]
	where the inverse exists. The matrix-valued Paley--Wiener theorem for convolution resolvents states that the resolvent belongs to $L^1$ precisely when $I_{S_d}-\widehat K(z)$ is invertible throughout the closed right half-plane. Since $S_d$ is finite-dimensional, invertibility is equivalent to non-vanishing of the determinant.
\end{proof}

We record two convolution estimates.
\begin{lemma}\label{lem:conv1}
	Let $m\in\mathbb N_0$, $z\in\C$, and let $C$ be a constant matrix. Suppose $g$ is exponentially bounded with $\|g(t)\|\leq Me^{\gamma t}$ for some $\gamma<\Re z$. Then
	\[
	\int_0^t(t-s)^m e^{z(t-s)}Cg(s)\,ds
	=t^me^{zt}C\widehat g(z)+o(t^me^{(\Re z)t}).
	\]
\end{lemma}
\begin{proof}
	Divide the integral by $t^me^{zt}$. The resulting integrand is $C(1-s/t)^me^{-zs}g(s)$ on $[0,t]$, which converges pointwise to $Ce^{-zs}g(s)$ and is dominated on $[0,\infty)$ by $\|C\|Me^{-(\Re z-\gamma)s}$. Dominated convergence proves the assertion.
\end{proof}

\begin{lemma}\label{lem:conv2}
	Let $m\in\mathbb N_0$, $z\in\C$, and $\alpha=\Re z$. Suppose
	\[
	f(t)=t^me^{zt}c+o(t^me^{\alpha t})
	\]
	for a constant vector $c$. If $\|R(t)\|\leq Me^{\alpha_Rt}$ for some $\alpha_R<\alpha$, then
	\[
	(R*f)(t)=t^me^{zt}\widehat R(z)c+o(t^me^{\alpha t}).
	\]
\end{lemma}
\begin{proof}
	Apply Lemma~\ref{lem:conv1} to the leading term. For the remainder $\varepsilon(t)=o(t^me^{\alpha t})$, split the convolution at a sufficiently large fixed $T$. The part over $[0,T]$ is bounded by $\int_0^T\|R(t-s)\|\,\|\varepsilon(s)\|\,ds=O(e^{\alpha_Rt})$, hence negligible. Given $\eta>0$, choose $T$ so that $\|\varepsilon(s)\|\leq\eta s^me^{\alpha s}$ for $s\geq T$. The remaining part, divided by $t^me^{\alpha t}$, is bounded by $M\eta\int_T^t e^{-(\alpha-\alpha_R)(t-s)}(s/t)^m\,ds\leq \frac{M\eta}{\alpha-\alpha_R}$. Letting $\eta\downarrow0$ proves the result.
\end{proof}

\section{Characteristic equation and dominant asymptotics}
Define
\[
x(t)=\E[X(t)\otimes X(t)],\qquad y(t)=\E[Y(t)\otimes Y(t)],
\]
and
\[
v_0(t)=y_0(t)\otimes y_0(t),\qquad \Gamma(t)=G(t)\otimes G(t).
\]
For every random vector $U\in\R^d$ with finite second moments, $\E[U\otimes U]\in S_d$ because it is a limit of finite sums of tensors of the form $u\otimes u$ (corresponding to a symmetric positive-semidefinite matrix). Hence both $x(t)$ and $y(t)$ take values in $S_d$.

\begin{theorem}\label{thm:Volterra}
	For every $t\geq0$,
	\begin{align}
		x(t)&=x_0(t)\otimes x_0(t)+(R*y)(t),\label{eq:xVolterra}\\
		y(t)&=v_0(t)+(\Gamma*y)(t).\label{eq:yVolterra}
	\end{align}
\end{theorem}
\begin{proof}
	For a deterministic matrix-valued function $K=(k_{ij})$ and a square-integrable adapted process $Z=(Z_j)$, the $(i,\ell)$ entry of 
	\[
	\E\bigl[\bigl(\int_0^t K(s)Z(s)\,dB(s)\bigr)\otimes\bigl(\int_0^t K(s)Z(s)\,dB(s)\bigr)\bigr]
	\] is
	\begin{multline*}\sum_{j,q}\E\left[ 	\left(\int_0^t k_{ij}(s)Z_j(s)\,dB(s)\right) 	\left(\int_0^t k_{\ell q}(s)Z_q(s)\,dB(s)\right)\right] 	\\=\int_0^t\sum_{j,q}k_{ij}(s)k_{\ell q}(s)\E[Z_j(s)Z_q(s)]\,ds
	\end{multline*}
	by the scalar It\^o isometry and bilinearity, which is the $(i,\ell)$ entry of $\int_0^t(K(s)\otimes K(s))\E[Z(s)\otimes Z(s)]\,ds$.
	
	Write $X(t)=x_0(t)+I_X(t)$, where $I_X(t)=\int_0^t r(t-s)Y(s)\,dB(s)$. The cross terms vanish after expectation because the stochastic integral has zero mean. Applying the tensor It\^o isometry with $K(s)=r(t-s)$ and $Z(s)=Y(s)$ gives
	\[\E[I_X(t)\otimes I_X(t)]=\int_0^t(r(t-s)\otimes r(t-s))y(s)\,ds=(R*y)(t),\]
	proving \eqref{eq:xVolterra}. Applying the same argument to \eqref{eq:Yvariation} with $K(s)=G(t-s)$ proves \eqref{eq:yVolterra}.
\end{proof}

The characteristic matrix of \eqref{eq:yVolterra}, on $S_d$, is $I_{S_d}-\widehat\Gamma_s(z)$. We now compute its transform. Set
\begin{equation}\label{eq:Theta}
	\Theta=\Sigma_1e^{A\tau}\otimes\Sigma_2+\Sigma_2\otimes\Sigma_1e^{A\tau}+\Sigma_2\otimes\Sigma_2,
	\qquad D=A_2+\Sigma_1\otimes\Sigma_1.
\end{equation}

\begin{theorem}\label{thm:transform}
	For $\Re z>\alpha_A$,
	\begin{equation}\label{eq:Gammahat}
		\widehat\Gamma(z)=\bigl(\Sigma_1\otimes\Sigma_1+e^{-z\tau}\Theta\bigr)(zI_{d^2}-A_2)^{-1}.
	\end{equation}
	Consequently, for $z\notin\spec(A_2)$, the meromorphic continuation satisfies
	\begin{equation}\label{eq:detidentity}
		\det\bigl(I_{d^2}-\widehat\Gamma(z)\bigr)
		=\frac{\det\bigl(zI_{d^2}-D-e^{-z\tau}\Theta\bigr)}{\det(zI_{d^2}-A_2)}.
	\end{equation}
	The same identities hold after restriction to $S_d$. In particular, in the half-plane $\Re z>\alpha_A$, the physical characteristic roots
	\[
	\mathcal Z_s=\{z\in\C:\Re z>\alpha_A,\ c_s(z)=0\},\qquad
	c_s(z)=\det\bigl(I_{S_d}-\widehat\Gamma_s(z)\bigr),
	\]
	are exactly the characteristic roots of the retarded equation
	\begin{equation}\label{eq:retarded}
		w'(t)=D_sw(t)+\Theta_sw(t-\tau)
	\end{equation}
	on $S_d$.
\end{theorem}
\begin{proof}
	Expand $G(t)\otimes G(t)$ into its four terms. For $\Re z>\alpha_A$, the instantaneous term has transform $(\Sigma_1\otimes\Sigma_1)(zI_{d^2}-A_2)^{-1}$. For the first cross term, the substitution $s=t-\tau$ and the identity $e^{As}\otimes e^{As}=e^{A_2s}$ yield
	\begin{align*} 	\int_\tau^\infty e^{-zt}(\Sigma_1e^{At})\otimes(\Sigma_2e^{A(t-\tau)})\,dt 	&=e^{-z\tau}(\Sigma_1e^{A\tau}\otimes\Sigma_2)\int_0^\infty e^{-zs}e^{A_2s}\,ds\\ 	&=e^{-z\tau}(\Sigma_1e^{A\tau}\otimes\Sigma_2)(zI_{d^2}-A_2)^{-1}, 	\end{align*}
	where the resolvent remains on the right because the coefficient matrices need not commute. The other cross term and the delayed-delayed term give the remaining two terms in $e^{-z\tau}\Theta(zI_{d^2}-A_2)^{-1}$. This proves \eqref{eq:Gammahat}.
	
	For $z\notin\spec(A_2)$, factor
	\[I_{d^2}-\widehat\Gamma(z)= 	\bigl[zI_{d^2}-A_2-\Sigma_1\otimes\Sigma_1-e^{-z\tau}\Theta\bigr](zI_{d^2}-A_2)^{-1},\]
	which gives \eqref{eq:detidentity}. The characteristic matrix of \eqref{eq:retarded} is $zI_{S_d}-D_s-e^{-z\tau}\Theta_s$. Since $\Re z>\alpha_A$ excludes $\spec(A_{2,s})$, the determinant identity gives equality of the two root sets in this half-plane, including multiplicities as zeros of the corresponding analytic determinants.
	
	It remains to check invariance of $S_d$. For any $u\in\R^d$, differentiating $e^{A_2t}(u\otimes u)=(e^{At}u)\otimes(e^{At}u)\in S_d$ at $t=0$ shows that $A_2(u\otimes u)=(Au)\otimes u+u\otimes(Au)\in S_d$. Likewise $\Sigma_1\otimes\Sigma_1$ and $\Sigma_2\otimes\Sigma_2$ preserve $S_d$. For the cross terms, if $a=\Sigma_1e^{A\tau}u$ and $b=\Sigma_2u$, then
	\[a\otimes b+b\otimes a=(a+b)\otimes(a+b)-a\otimes a-b\otimes b\in S_d.\]
	Thus $\Theta$, $D$, and $A_2$ preserve $S_d$, and the same calculation applies to their restrictions.
\end{proof}

\begin{remark}[The exceptional set]\label{rem:exceptional}
	The quotient identity \eqref{eq:detidentity} is asserted only away from $\spec(A_2)$. At a point of $\spec(A_2)$, zeros of the numerator of the quotient need not correspond to zeros of the meromorphic continuation of the left-hand determinant: cancellation between numerator and denominator may occur. This is why the physical characteristic roots are defined in the half-plane $\Re z>\alpha_A$, where the Laplace transform is originally defined and the denominator is nonzero. All roots used in the dominant asymptotic theorem and in the stability criterion lie in this half-plane. In particular, once $A$ is Hurwitz, the closed right half-plane is contained in it.
\end{remark}

Let
\[
\alpha_d=\sup\{\Re z:z\in\mathcal Z_s\},\qquad \sup\varnothing=-\infty.
\]
Thus $\alpha_A$ is the spectral bound of the instantaneous tensor generator, while $\alpha_d$ is the spectral bound of the physical feedback characteristic roots. The roots of a retarded characteristic equation are isolated and only finitely many lie in any right half-plane; see Hale and Verduyn Lunel \cite[Chapter~8]{Hale} and Diekmann et al.\ \cite{Diekmann}. Consequently, if $\alpha_d>\alpha_A$, then $\alpha_d$ is attained by finitely many physical roots.

\begin{theorem}[Dominant mean-square asymptotics]\label{thm:asymptotic}
	Suppose $\alpha_d>\alpha_A$, and let $z_1,\ldots,z_K$ be the physical characteristic roots with real part $\alpha_d$. Let $m_k$ be the pole order of $(I_{S_d}-\widehat\Gamma_s(z))^{-1}$ at $z_k$, and put
	\[
	m_* =\max_{1\leq k\leq K}m_k,\qquad K_* =\{k:m_k=m_*\}.
	\]
	For $k\in K_*$ define
	\[
	C_{k,m_*-1}=\frac1{(m_*-1)!}\lim_{z\to z_k}(z-z_k)^{m_*}
	(I_{S_d}-\widehat\Gamma_s(z))^{-1}.
	\]
	Then, as $t\to\infty$,
	\begin{equation}\label{eq:dominant}
		x(t)=t^{m_*-1}\sum_{k\in K_*}e^{z_kt}(z_kI_{S_d}-A_{2,s})^{-1}
		C_{k,m_*-1}\widehat v_0(z_k)+o(t^{m_*-1}e^{\alpha_dt}).
	\end{equation}
\end{theorem}
\begin{proof}
	Let $\rho_s$ be the Volterra resolvent of $\Gamma_s$. Its transform satisfies
	\[I_{S_d}+\widehat\rho_s(z)=(I_{S_d}-\widehat\Gamma_s(z))^{-1}.\]
	Hence the physical characteristic roots are precisely the poles of the transformed resolvent in $\Re z>\alpha_A$. Let $H_s$ be the fundamental matrix of \eqref{eq:retarded}, with $H_s(t)=0$ for $t<0$ and $H_s(0)=I_{S_d}$. Its Laplace transform is
	\[\widehat H_s(z)=(zI_{S_d}-D_s-e^{-z\tau}\Theta_s)^{-1}.\]
	Writing $B_s=(\Sigma_1\otimes\Sigma_1)_s=D_s-A_{2,s}$ and using \eqref{eq:Gammahat}, we have
	\[\widehat\rho_s(z)=(zI_{S_d}-A_{2,s})\widehat H_s(z)-I_{S_d}=B_s\widehat H_s(z)+e^{-z\tau}\Theta_s\widehat H_s(z).\]
	Inverting the transform yields the identity
	\begin{equation}\label{eq:rhoH}
		\rho_s(t)=B_sH_s(t)+\Theta_sH_s(t-\tau),\qquad t\geq0.
	\end{equation}
	
	Choose $\beta\in(\alpha_A,\alpha_d)$ so that the only characteristic roots with real part greater than $\beta$ are $z_1,\ldots,z_K$. The standard spectral expansion for retarded equations (see, e.g., Diekmann et al.\ \cite[Chapter~I, Theorem~5.4]{Diekmann} and Hale and Verduyn Lunel \cite[Chapter~7]{Hale}) gives
	\[H_s(t)=\sum_{k=1}^K e^{z_kt}\sum_{j=0}^{m_k-1}t^jH_{k,j}+O(e^{\beta t}),\]
	where the coefficients are the corresponding principal-part coefficients. By \eqref{eq:rhoH}, the same expansion holds for $\rho_s$, with an exponentially smaller remainder. In particular, the coefficient of $t^{m_k-1}e^{z_kt}$ is $C_{k,m_k-1}$ as defined above. Thus
	\[\rho_s(t)=\sum_{k=1}^K e^{z_kt}\sum_{j=0}^{m_k-1}t^jC_{k,j}+O(e^{\beta t}).\]
	Here and below the $C_{k,j}$ for lower powers denote the corresponding coefficients, without requiring an explicit formula for them.
	
	Since $y_0(t)$ is a linear combination of $e^{At}\psi(0)$ and the initial history on $[0,\tau]$, there is, for every $\eta>0$, a constant $C_\eta$ (depending on $\psi$) such that
	\[\|v_0(t)\|+\|R(t)\|\leq C_\eta e^{(\alpha_A+\eta)t}.\]
	Choose $\eta>0$ with $\alpha_A+\eta<\alpha_d$. Since $y=v_0+\rho_s*v_0$, Lemma~\ref{lem:conv1} applied to the dominant terms gives
	\[ 	y(t)=t^{m_*-1}\sum_{k\in K_*}e^{z_kt}C_{k,m_*-1}\widehat v_0(z_k) 	+o(t^{m_*-1}e^{\alpha_dt}). 	\]
	Terms of lower polynomial degree and the remainder with exponential rate below $\alpha_d$ are absorbed into the stated error. Applying Lemma~\ref{lem:conv2} to $x=x_0\otimes x_0+R*y$, and using \eqref{eq:Rhat}, gives \eqref{eq:dominant}. The term $x_0(t)\otimes x_0(t)$ has exponential order at most $\alpha_A+\eta<\alpha_d$ and is negligible.
\end{proof}

\begin{remark}[Pole order and characteristic multiplicity]\label{rem:pole}
	The integer $m_k$ in Theorem~\ref{thm:asymptotic} is the pole order of the matrix resolvent $(I_{S_d}-\widehat\Gamma_s(z))^{-1}$. It is this pole order that determines the polynomial factor $t^{m_k-1}$ in the time-domain expansion. It should not be identified without qualification with the algebraic multiplicity of $z_k$ as a zero of a characteristic determinant, or with the multiplicity of a characteristic root in a particular scalar factorisation. The resolvent pole order is the invariant quantity relevant to the asymptotic expansion.
\end{remark}

\section{Accessibility, generic non-cancellation and stability}
A characteristic root with nonnegative real part prevents mean-square stability only if its associated dominant mode can be excited by an admissible initial history. We make this issue explicit. For a subspace $U\subseteq\R^d$, put
\[
\Sym(U)=\operatorname{span}\{u\otimes u:u\in U\}\subseteq S_d.
\]

\begin{lemma}[Accessibility of a dominant mode]\label{lem:access}
	Suppose $\alpha_d>\alpha_A$, and let $z_k$ be a dominant physical root with pole order $m_*$. If
	\begin{equation}\label{eq:access}
		C_{k,m_*-1}\Sym(\Ran\Sigma_2)\ne\{0\},
	\end{equation}
	then there is an admissible continuous initial history for which the coefficient of the $z_k$-mode in \eqref{eq:dominant} is nonzero. In particular, \eqref{eq:access} holds whenever $\Sigma_2$ is nonsingular.
\end{lemma}
\begin{proof}
	Choose $u\in\Ran\Sigma_2$ such that $C_{k,m_*-1}(u\otimes u)\ne0$, and choose $w\in\R^d$ with $\Sigma_2w=u$. Fix $t_0\in(0,\tau)$. By continuity, there is a compact subinterval $I\subset(0,\tau)$ containing $t_0$ in its interior on which
	\[|e^{-z_kt}-e^{-z_kt_0}|<\tfrac12|e^{-z_kt_0}|.\]
	Choose a nonzero continuous real-valued function $g$ supported in $I$, and define
	\[\psi(\theta)=g(\theta+\tau)w\quad(-\tau\leq\theta<0),\qquad \psi(0)=0.\]
	Then $x_0(t)=0$ for $t\geq0$, and $y_0(t)=g(t)u$ on $[0,\tau]$, with $y_0(t)=0$ for $t>\tau$. Therefore
	\[\widehat v_0(z_k)=\left(\int_0^\tau e^{-z_kt}g(t)^2\,dt\right)(u\otimes u).\]
	The choice of $I$ ensures that the scalar integral is nonzero: its difference from $e^{-z_kt_0}\int_0^\tau g(t)^2\,dt$ has modulus less than half the modulus of that quantity. Hence $C_{k,m_*-1}\widehat v_0(z_k)\ne0$. Since $\Re z_k=\alpha_d>\alpha_A$, the matrix $z_kI_{S_d}-A_{2,s}$ is invertible, so the coefficient in \eqref{eq:dominant} is nonzero. Finally, if $\Sigma_2$ is nonsingular, then $\Sym(\Ran\Sigma_2)=S_d$; since $S_d^{\C}=\operatorname{span}_{\C}(S_d)$ and $C_{k,m_*-1}\ne0$ on $S_d^{\C}$, it cannot annihilate all of $S_d$.
\end{proof}

\begin{theorem}[Generic non-cancellation]\label{thm:generic}
	Suppose $\alpha_d>\alpha_A$. Define the leading-coefficient map
	\[
	Q(\psi)=\bigl(C_{k,m_*-1}\widehat v_{0,\psi}(z_k)\bigr)_{k\in K_*}
	\in\prod_{k\in K_*}S_d^{\C}.
	\]
	The actual coefficients in the expansion \eqref{eq:dominant} are obtained by applying the invertible factors $(z_kI_{S_d}-A_{2,s})^{-1}$ componentwise, so they have the same common zero set as $Q$. If $Q$ is not identically zero on $\cC$, then
	\[
	\{\psi\in\cC:Q(\psi)=0\}
	\]
	is closed and has empty interior. Thus the histories for which at least one dominant leading coefficient is nonzero form an open dense subset of $\cC$.
\end{theorem}
\begin{proof}
	The map $\psi\mapsto y_{0,\psi}$ is linear and continuous. The map
	\[V(\psi,\phi)(t)=\tfrac12\bigl(y_{0,\psi}(t)\otimes y_{0,\phi}(t)+y_{0,\phi}(t)\otimes y_{0,\psi}(t)\bigr)\]
	is continuous, symmetric, and bilinear, and $V(\psi,\psi)=v_{0,\psi}$. For each dominant root $z_k$, the Laplace functional $v\mapsto\widehat v(z_k)$ is continuous on these functions because $y_0$ grows at most like $e^{(\alpha_A+\eta)t}$ and $\Re z_k=\alpha_d>\alpha_A$. It follows that $Q$ is a continuous homogeneous quadratic map from $\cC$ into a finite-dimensional complex vector space. Its zero set is therefore closed.
	
	If the zero set contained an open ball centred at $\psi_0$, choose a complex linear functional $\ell$ such that $q=\ell\circ Q$ is not identically zero. For any $h\in\cC$, the scalar polynomial $t\mapsto q(\psi_0+th)$, for real $t$ sufficiently close to zero, would vanish on an interval. Its coefficient of $t^2$ is $q(h)$, so $q(h)=0$. Since $h$ is arbitrary, this contradicts $q\not\equiv0$. The zero set has empty interior.
\end{proof}

\begin{theorem}[Mean-square asymptotic stability]\label{thm:stability}
	Assume the following accessibility condition: whenever $\mathcal Z_s$ contains a root in the closed right half-plane, at least one root $z_k$ with $\Re z_k=\alpha_d$ and pole order $m_*$ satisfies \eqref{eq:access}. Under this condition, the zero solution of \eqref{eq:model} is mean-square asymptotically stable for every deterministic initial history $\psi\in\cC$, in the sense that
	\[
	\E\|X(t)\|^2\longrightarrow0\qquad(t\to\infty),
	\]
	if and only if
	\begin{enumerate}
		\item[(i)] $A$ is Hurwitz, and
		\item[(ii)] $c_s(z)\ne0$ for every $z$ with $\Re z\geq0$.
	\end{enumerate}
	Moreover, when \textup{(i)} and \textup{(ii)} hold, the uniform bound \eqref{eq:uniform_bound} yields mean-square stability in addition to attraction. In particular, the accessibility condition is automatic if $\Sigma_2$ is nonsingular.
\end{theorem}
\begin{proof}
	Suppose first that (i) and (ii) hold. Since $A$ is Hurwitz, $A_2$ is Hurwitz and $\alpha_A<0$. Also $G$ decays exponentially, so $\Gamma_s\in L^1(0,\infty;\cL(S_d))$. Condition (ii) and Lemma~\ref{lem:PW} imply that its Volterra resolvent $\rho_s$ is integrable. Setting $\gamma=-\alpha_A/4>0$, there exist constants $C_0,C_1>0$ independent of $\psi$ such that $\|x_0(t)\|\leq C_0e^{-\gamma t}\|\psi\|_{\cC}$ and $\|y_0(t)\|\leq C_1e^{-\gamma t}\|\psi\|_{\cC}$ for all $t\geq0$, whence
	\[\|x_0(t)\otimes x_0(t)\|+\|v_0(t)\|\leq C_2e^{-2\gamma t}\|\psi\|_{\cC}^2,\qquad t\geq0.\]
	Since $y=v_0+\rho_s*v_0$ and $\rho_s\in L^1(0,\infty;\cL(S_d))$, we obtain $\|y(t)\|\leq \eta_1(t)\|\psi\|_{\cC}^2$, where $\eta_1(t)=C_2\bigl(e^{-2\gamma t}+\int_0^t\|\rho_s(t-s)\|e^{-2\gamma s}\,ds\bigr)$ is bounded on $[0,\infty)$ and satisfies $\eta_1(t)\to0$ as $t\to\infty$. Applying the same estimate to $x=x_0\otimes x_0+R_s*y$, where $\|R_s(t)\|\leq M_Re^{-2\gamma t}$, gives
	\begin{equation}\label{eq:uniform_bound}
		\E\|X(t)\|^2\leq \eta(t)\|\psi\|_{\cC}^2,\qquad t\geq0,
	\end{equation}
	for a deterministic function $\eta:[0,\infty)\to[0,\infty)$ independent of $\psi$ with $\sup_{t\geq0}\eta(t)<\infty$ and $\eta(t)\to0$ as $t\to\infty$. In particular, $\E\|X(t)\|^2\to0$.
	
	Conversely, suppose $x(t)\to0$ for every deterministic initial history. Taking expectations in \eqref{eq:Xvariation} gives
	\[\E X(t)=e^{At}\psi(0).\]
	By Jensen's inequality, $\|\E X(t)\|^2\leq\E\|X(t)\|^2\to0$. Since $\psi(0)$ is arbitrary, $e^{At}\xi\to0$ for every $\xi\in\R^d$, and therefore $A$ is Hurwitz. In particular, $\alpha_A<0$.
	
	If (ii) failed, there would be a physical characteristic root in the closed right half-plane, so $\alpha_d\geq0>\alpha_A$. Theorem~\ref{thm:asymptotic} applies. By the accessibility hypothesis and Lemma~\ref{lem:access}, the leading coefficient map in Theorem~\ref{thm:generic} is not identically zero. By Theorem~\ref{thm:generic}, choose an initial history for which at least one dominant leading coefficient is nonzero. Combining terms with equal imaginary parts, the leading term in \eqref{eq:dominant} has the form
	\[t^{m_*-1}e^{\alpha_dt}F(t),\qquad F(t)=\sum_{\ell=1}^N e^{i\omega_\ell t}M_\ell,\]
	where the frequencies $\omega_\ell$ are distinct and not all $M_\ell$ vanish. In a Euclidean inner product $\langle\cdot,\cdot\rangle$ on $S_d$, orthogonality of distinct exponentials under Ces\`aro averaging gives
	\[\lim_{T\to\infty}\frac1T\int_0^T\|F(t)\|^2\,dt 	=\sum_{\ell,q=1}^N\left(\lim_{T\to\infty}\frac1T\int_0^T e^{i(\omega_\ell-\omega_q)t}\,dt\right)\langle M_\ell,M_q\rangle 	=\sum_{\ell=1}^N\|M_\ell\|^2>0.\]
	Thus there are $c>0$ and $t_n\to\infty$ such that $\|F(t_n)\|\geq c$. Since $\alpha_d\geq0$ and $m_*-1\geq0$, the asymptotic expansion implies that $\|x(t_n)\|$ cannot tend to zero. This contradicts the assumed mean-square convergence for every history. Therefore (ii) holds.
	
	If $\Sigma_2$ is nonsingular, Lemma~\ref{lem:access} supplies the required accessibility condition automatically.
\end{proof}

\begin{corollary}[Random initial histories]\label{cor:random}
	Suppose that conditions \textup{(i)} and \textup{(ii)} of Theorem~\ref{thm:stability} hold, and let $\psi$ be a $\cC$-valued random initial history independent of the Brownian motion $B$ with $\E\|\psi\|_{\cC}^2<\infty$. Then the corresponding solution of \eqref{eq:model} satisfies
	\[
	\E\|X(t)\|^2\longrightarrow0\qquad(t\to\infty).
	\]
\end{corollary}
\begin{proof}
	Conditioning on $\psi$ and applying the uniform bound \eqref{eq:uniform_bound} pathwise in $\psi$ yields
	\[\E\|X(t)\|^2=\E\bigl[\E[\|X(t)\|^2\mid\psi]\bigr]\leq \eta(t)\,\E\|\psi\|_{\cC}^2,\qquad t\geq0.\]
	Since $\E\|\psi\|_{\cC}^2<\infty$ and $\eta(t)\to0$ as $t\to\infty$, the conclusion follows.
\end{proof}

\section{Examples and concluding remarks}

We conclude with two two-dimensional examples illustrating the spectral criterion when $\Sigma_2$ is nonsingular and the role of the delayed-channel accessibility condition \eqref{eq:access} when $\Sigma_2$ is singular.

\begin{example}[Multiple dominant roots with nonsingular $\Sigma_2$]\label{ex:nonsingular}
	Take $d=2$, $\tau=1$, and
	\[
	A=\begin{pmatrix}-1&-\pi\\ \pi&-1\end{pmatrix},\qquad
	\Sigma_1=0,\qquad \Sigma_2=\sqrt2 I_2.
	\]
	Then $\spec(A)=\{-1\pm\pi i\}$, so $A$ is Hurwitz. Moreover $D=A_2$ and $\Theta=2I_4$. On $S_2$, the eigenvalues of $A_{2,s}$ are $-2$ and $-2\pm2\pi i$ (the fourth eigenvalue $-2$ of $A_2$ lies on the antisymmetric subspace). The physical characteristic determinant therefore factors as
	\[
	\det\bigl(zI_{S_2}-A_{2,s}-2e^{-z}I_{S_2}\bigr)
	=(z+2-2e^{-z})(z+2-2\pi i-2e^{-z})(z+2+2\pi i-2e^{-z}).
	\]
	Thus $z=0,2\pi i,-2\pi i$ are characteristic roots. They are the only roots in the closed right half-plane. Indeed, each factor has the form $z+2-i\eta=2e^{-z}$ with $\eta\in\{0,2\pi,-2\pi\}$. If $z=x+iy$ and $x\geq0$, then
	\[
	2e^{-x}=|z+2-i\eta|=\sqrt{(x+2)^2+(y-\eta)^2}\geq x+2\geq2.
	\]
	Since $2e^{-x}\leq2$, equality holds throughout, forcing $x=0$ and $y=\eta$. Hence $\alpha_d=0$ and $\alpha_A=-2$.
	
	For $f_\eta(z)=z+2-i\eta-2e^{-z}$ one has $f_\eta'(i\eta)=3\ne0$, and since the three eigenvalues of $A_{2,s}$ are distinct, the three scalar factors do not vanish simultaneously at $z\in\{0,\pm2\pi i\}$. Hence all three roots correspond to simple poles ($m_*=1$) of $(I_{S_2}-\widehat\Gamma_s(z))^{-1}$. Theorem~\ref{thm:asymptotic} consequently gives
	\[
	x(t)=M_0+e^{2\pi it}M_++e^{-2\pi it}M_-+o(1),
	\]
	where the matrices depend on the initial history. Since $\Sigma_2$ is nonsingular, the accessibility condition holds automatically, and the presence of roots on the imaginary axis is consistent with the failure of mean-square asymptotic stability.
\end{example}

\begin{example}[Failure of delayed-channel accessibility for singular $\Sigma_2$]\label{ex:singular}
	Take $d=2$, $\tau=1$, and
	\[
	A=-I_2,\qquad
	\Sigma_1=\begin{pmatrix}0&0\\0&\sqrt2\end{pmatrix},\qquad
	\Sigma_2=\begin{pmatrix}1&0\\0&0\end{pmatrix}.
	\]
	Here $A$ is Hurwitz with $A_{2,s}=-2I_{S_2}$ and $\alpha_A=-2$. In the orthonormal basis $\bigl\{e_1\otimes e_1,\,\frac{1}{\sqrt2}(e_1\otimes e_2+e_2\otimes e_1),\,e_2\otimes e_2\bigr\}$ of $S_2$, we have
	\[
	D_s=\operatorname{diag}(-2,-2,0),\qquad
	\Theta_s=\operatorname{diag}(1,\sqrt2e^{-1},0).
	\]
	For $\Re z>-2$, the physical characteristic determinant is therefore
	\[
	c_s(z)=\det\bigl(I_{S_2}-\widehat\Gamma_s(z)\bigr)
	=\frac{(z+2-e^{-z})(z+2-\sqrt2e^{-1-z})z}{(z+2)^3}.
	\]
	In the closed right half-plane $\Re z\geq0$, the first two factors are non-vanishing because $|z+2|\geq2>1\geq|e^{-z}|$ and $|z+2|\geq2>\sqrt2e^{-1}\geq|\sqrt2e^{-1-z}|$. Thus $z_1=0$ is the unique root in $\{\Re z\geq0\}$, with $\alpha_d=0>-2=\alpha_A$ and pole order $m_*=1$. Near $z=0$,
	\[
	\bigl(I_{S_2}-\widehat\Gamma_s(z)\bigr)^{-1}
	=\operatorname{diag}\left(\frac{z+2}{z+2-e^{-z}},\,\frac{z+2}{z+2-\sqrt2e^{-1-z}},\,\frac{z+2}{z}\right),
	\]
	so $C_{1,0}=\lim_{z\to0}z(I_{S_2}-\widehat\Gamma_s(z))^{-1}=\operatorname{diag}(0,0,2)$. Since $\Ran\Sigma_2=\operatorname{span}\{e_1\}$, we have $\Sym(\Ran\Sigma_2)=\operatorname{span}\{e_1\otimes e_1\}$, and hence
	\[
	C_{1,0}\Sym(\Ran\Sigma_2)=\{0\}.
	\]
	Thus \eqref{eq:access} fails: the root $z_1=0$ cannot be excited through the delayed history channel $\Sym(\Ran\Sigma_2)$ alone as in Lemma~\ref{lem:access}. Nevertheless, any history with $\psi_2(0)\ne0$ excites the $e_2\otimes e_2$ mode via $\Sigma_1$ (giving $\E[X_2(t)^2]=\psi_2(0)^2\not\to0$), so mean-square asymptotic stability still fails. In general, a full necessary-and-sufficient accessibility condition for singular $\Sigma_2$ requires $\operatorname{span}\{\widehat v_{0,\psi}(z_k):\psi\in\cC\}\not\subseteq\ker C_{k,m_*-1}$, whereas \eqref{eq:access} provides an explicit, directly checkable sufficient criterion purely in terms of $\Sigma_2$.
\end{example}

\begin{remark}[Scope and multi-dimensional noise]\label{rem:scope}
	The results above treat the constant-coefficient model \eqref{eq:model} with instantaneous drift, a single discrete delay, and scalar Brownian motion. Replacing $B$ by an $m$-dimensional Brownian motion and the diffusion coefficient by $\sum_{j=1}^m(\Sigma_{1,j}X(t)+\Sigma_{2,j}X(t-\tau))\,dB_j(t)$ preserves the basic Volterra strategy: the kernel becomes $\Gamma(t)=\sum_{j=1}^mG_j(t)\otimes G_j(t)$, and the corresponding operators $D$ and $\Theta$ are obtained by summing the respective Kronecker terms. The accessibility condition must, however, be reformulated using the actual set of covariance profiles attainable from a common initial history; it does not follow without further argument by simply replacing $\Sym(\Ran\Sigma_2)$ with the sum of the individual symmetric range spaces. Extensions to delayed drift, multiple discrete delays, or distributed delay kernels require a different treatment of the deterministic variation-of-constants representation and lie beyond the scope of the present analysis.
\end{remark}


\begin{thebibliography}{99}
	\bibitem{ApplebyVolterra} J. A. D. Appleby, Mean square characterisation of a stochastic Volterra integrodifferential equation with delay, \emph{International Journal of Dynamical Systems and Differential Equations} \textbf{11} (2021), no. 3--4, 194--226.
	\bibitem{ApplebyLawless} J. A. D. Appleby and E. Lawless, Mean square asymptotic stability characterisation of perturbed linear stochastic functional differential equations, \emph{Applied Numerical Mathematics} \textbf{200} (2024), 80--109.
	\bibitem{AMR} J. A. D. Appleby, X. Mao and M. Riedle, Geometric Brownian motion with delay: Mean square characterisation, \emph{Proceedings of the American Mathematical Society} \textbf{137} (2009), no. 1, 339--348.
	\bibitem{Bao} J. Bao, G. Yin and C. Yuan, \emph{Asymptotic Analysis for Functional Stochastic Differential Equations}, SpringerBriefs in Mathematics, Springer, Cham, 2016.
	\bibitem{Bellman54} R. Bellman, Limit theorems for non-commutative operations. I, \emph{Duke Mathematical Journal} \textbf{21} (1954), no. 3, 491--500.
	\bibitem{Bellman59} R. Bellman, Kronecker products and the second method of Lyapunov, \emph{Mathematische Nachrichten} \textbf{20} (1959), no. 1--2, 17--19.
	\bibitem{BellmanCooke} R. Bellman and K. L. Cooke, \emph{Differential-Difference Equations}, Academic Press, New York, 1963.
	\bibitem{BuckwarNotarangelo} E. Buckwar and G. Notarangelo, A note on the analysis of asymptotic mean-square stability properties for systems of linear stochastic delay differential equations, \emph{Discrete and Continuous Dynamical Systems -- Series B} \textbf{18} (2013), no. 6, 1521--1531.
	\bibitem{BuckwarSickenberger} E. Buckwar and T. Sickenberger, A structural analysis of asymptotic mean-square stability for multi-dimensional linear stochastic differential systems, \emph{Applied Numerical Mathematics} \textbf{62} (2012), no. 7, 842--859.
	\bibitem{Diekmann} O. Diekmann, S. A. van Gils, S. M. Verduyn Lunel and H.-O. Walther, \emph{Delay Equations: Functional-, Complex-, and Nonlinear Analysis}, Applied Mathematical Sciences 110, Springer, New York, 1995.
	\bibitem{Gripenberg} G. Gripenberg, S.-O. Londen and O. Staffans, \emph{Volterra Integral and Functional Equations}, Encyclopedia of Mathematics and its Applications, Cambridge University Press, Cambridge, 1990.
	\bibitem{Hale} J. K. Hale and S. M. Verduyn Lunel, \emph{Introduction to Functional Differential Equations}, Applied Mathematical Sciences 99, Springer, New York, 1993.
	\bibitem{IklodiDankowicz} Z. Iklodi and H. Dankowicz, Algebraic conditions for second-moment stability boundaries of linear, time-invariant stochastic delay-differential equations, arXiv:2607.01374, 2026.
	\bibitem{Khasminskii} R. Z. Khasminskii, \emph{Stochastic Stability of Differential Equations}, Sijthoff \& Noordhoff, Alphen aan den Rijn, 1980.
	\bibitem{KolmanovskiiMyshkis} V. B. Kolmanovskii and A. Myshkis, \emph{Introduction to the Theory and Applications of Functional-Differential Equations}, Kluwer Academic Publishers, Dordrecht, 1999.
	\bibitem{KolmanovskiiNosov} V. B. Kolmanovskii and V. R. Nosov, \emph{Stability of Functional-Differential Equations}, Academic Press, London, 1986.
	\bibitem{LeiMackey} J. Lei and M. C. Mackey, Stochastic differential delay equation, moment stability, and application to hematopoietic stem cell regulation system, \emph{SIAM Journal on Applied Mathematics} \textbf{67} (2007), no. 2, 387--407.
	\bibitem{MaoBook} X. Mao, \emph{Stochastic Differential Equations and Applications}, 2nd ed., Horwood Publishing, Chichester, 2007.
	\bibitem{MaoRaz} X. Mao, Razumikhin-type theorems on exponential stability of stochastic functional differential equations, \emph{Stochastic Processes and their Applications} \textbf{65} (1996), no. 2, 233--250.
	\bibitem{MaoShah} X. Mao and A. Shah, Exponential stability of stochastic differential delay equations, \emph{Stochastics and Stochastics Reports} \textbf{60} (1997), no. 1--2, 135--153.
	\bibitem{MaoRiedle} X. Mao and M. Riedle, Mean square stability of stochastic Volterra integro-differential equations, \emph{Systems \& Control Letters} \textbf{55} (2006), no. 6, 459--465.
	\bibitem{MackeyNechaeva} M. C. Mackey and I. G. Nechaeva, Solution moment stability in stochastic differential delay equations, \emph{Physical Review E} \textbf{52} (1995), no. 4, 3366--3376.
	\bibitem{Mohammed} S. E. A. Mohammed, \emph{Stochastic Functional Differential Equations}, Research Notes in Mathematics 99, Pitman, London, 1984.
	\bibitem{MizelTrutzer} V. J. Mizel and V. Trutzer, Stochastic heredity equations: existence and asymptotic stability, \emph{Journal of Integral Equations} \textbf{7} (1984), 1--72.
	\bibitem{NevelsonKhasminskii} M. B. Nevelson and R. Z. Khasminskii, Stability of a linear system with random disturbances of its parameters, \emph{Journal of Applied Mathematics and Mechanics} \textbf{30} (1966), no. 2, 487--493.
	\bibitem{Shaikhet} L. Shaikhet, \emph{Lyapunov Functionals and Stability of Stochastic Functional Differential Equations}, Springer, Cham, 2013.
	\bibitem{WangLiLei} Z. Wang, X. Li and J. Lei, Second moment boundedness of linear stochastic delay differential equations, \emph{Discrete and Continuous Dynamical Systems -- Series B} \textbf{19} (2014), no. 9, 2963--2991.
	\bibitem{WangDistributed} Z. Wang, X. Li and J. Lei, Moment boundedness of linear stochastic delay differential equations with distributed delay, \emph{Stochastic Processes and their Applications} \textbf{124} (2014), no. 1, 586--612.
\end{thebibliography}
\end{document}